\documentstyle{amltd2004}
\begin{document}
\input boxedeps.tex 
\SetepsfEPSFSpecial 
\HideDisplacementBoxes
\def\figin#1#2{
$$
 {\BoxedEPSF{#1.eps scaled
#2}%
}%
$$
\noindent}
\annalsline{155}{2002}
\received{November 2, 2000}
\startingpage{895}

\def\joinrel{\mathrel{\mkern-4mu}}
\def\relbar{\mathrel{\smash-}}
\def\lrar{\relbar\joinrel\relbar\joinrel\rightarrow}
\def\llar{\leftarrow\joinrel\relbar\joinrel\relbar}

\def\bye{\end{document}}
 \font\tenrm=cmr10
\input amssym.def
\input amssym.tex
\catcode`\@=11
\font\twelvemsb=msbm10 scaled 1100
\font\tenmsb=msbm10
\font\ninemsb=msbm10 scaled 800
\newfam\msbfam
\textfont\msbfam=\twelvemsb  \scriptfont\msbfam=\ninemsb
  \scriptscriptfont\msbfam=\ninemsb
\def\msb@{\hexnumber@\msbfam}
\def\Bbb{\relax\ifmmode\let\next\Bbb@\else
 \def\next{\errmessage{Use \string\Bbb\space only in math
mode}}\fi\next}
\def\Bbb@#1{{\Bbb@@{#1}}}
\def\Bbb@@#1{\fam\msbfam#1}
\catcode`\@=12

 \catcode`\@=11
\font\twelveeuf=eufm10 scaled 1100
\font\teneuf=eufm10
\font\nineeuf=eufm7 scaled 1100
\newfam\euffam
\textfont\euffam=\twelveeuf  \scriptfont\euffam=\teneuf
  \scriptscriptfont\euffam=\nineeuf
\def\euf@{\hexnumber@\euffam}
\def\frak{\relax\ifmmode\let\next\frak@\else
 \def\next{\errmessage{Use \string\frak\space only in math
mode}}\fi\next}
\def\frak@#1{{\frak@@{#1}}}
\def\frak@@#1{\fam\euffam#1}
\catcode`\@=12

\newcommand{\bo}[1]{{\bf #1}}
\newcommand{\lra}{\longrightarrow}
\newcommand{\ra}{\rightarrow}
\newcommand{\Spaces}{\bo{Spaces}}
\newcommand{\algch}{\bo{Alg}^h_{\C}}
\newcommand{\algc}{\bo{Alg_C}}
\newcommand{\alg}{\bo{Alg_T}}
\newcommand{\T}{{\bo T}}
\newcommand{\C}{{\bo C}}
\newcommand{\M}{{\bo M}}
\newcommand{\gr}{\bo {Gr}}
\newcommand{\HoM}{{{\bo H}{\bo o}\M}}
\newcommand{\sC}{\Spaces^{\C}}
\newcommand{\st}{\Spaces^{\T}}
\newcommand{\stf}{\Spaces^{\T}_{\rm fib}}
\newcommand{\stc}{\Spaces^{\T}_{\rm cof}}
\newcommand{\scf}{\Spaces^{\C}_{\rm fib}}
\newcommand{\scc}{\Spaces^{\C}_{\rm cof}}
\newcommand{\Cdi}{\C^{\rm disc}}
\newcommand{\sCdi}{\Spaces^{\Cdi}}
\newcommand{\lst}{{{\bo L}\st}}
\def\mapright#1#2{\smash{\mathop{\hbox to
#1pt{\rightarrowfill}}\limits^{#2}}}
\newcommand{\Hom}{{\rm Hom}}
\newcommand{\Map}{{\rm RMap}}
\newcommand{\map}{{\rm Map}}

\title{Algebraic theories in homotopy theory} 
\shorttitle{Algebraic theories in homotopy theory} 

  \author{Bernard Badzioch}
 \institutions{University of Minnesota, 
Minneapolis, MN\\
{\eightpoint {\it E-mail address\/}: badzioch@math.umn.edu
}}

\section{Introduction}

It is well known in homotopy theory that 
given a loop space $ X$ one can always find 
a simplicial group $ G$ weakly equivalent 
to $ X$, such that the weak equivalence 
can be realized by maps
preserving multiplication. It is also known 
that loop spaces are not the only class of spaces 
for which a statement of this kind holds; for example, 
any $ A_{\infty}$-space is weakly equivalent to 
a simplicial monoid and every Eilenberg-Mac Lane
space $ K(G,n)$ with $G$  an abelian group is equivalent
to a simplicial abelian monoid. 
Results like this suggest that there might be
some general principle comparing 
homotopy structures  on a space
to algebraic structures. 
Our aim in this paper is to show that there is, in fact, 
such a principle. 
To make this precise we need a few definitions. 

\numbereddemo{Definition}
\label{STRICTALG}
An {\it algebraic theory} $ \T$ is  a {small category }
with objects 
$T_0, T_1, \dots$ together with, for each $n$, an expression of 
$T_n$ as the categorical product in $\T$ of $n$ copies of 
the object $T_1$. In particular $ T_0$ is the {terminal object}
in $ \T$. We assume that it is also  the initial object.
   
Given an algebraic theory $\T$, a {\it strict} $\T$-\/{\it algebra} $A$ is 
a {product-preserving} functor $A\colon\T\ra \Spaces$. 
\enddemo  

We will denote by $ \alg$ the category of all strict 
$ \T$-algebras with natural transformations of functors 
as morphisms. A strict
$\T$-algebra structure on a space $Y$ is a strict $\T$-algebra $A$
together with an isomorphism $Y\cong A(T_1)$.

Algebraic theories appear naturally in the
study of algebraic structures. For example, 
let $ \gr$ be the category of groups and for 
$ n\geq 0$ let $ F_n$ denote the free group generated 
by the set $ \{1,\dots,n\}$ ($F_0$ is the trivial group).
Define $ \T_\gr^{\rm op}$ to be the full subcategory 
of $ \gr$ with objects 
$ F_0,F_1,\dots\ $. Its opposite category $ \T_\gr$  
is then an algebraic theory.
To see this, observe that $n$ inclusions
$\{1\}\hookrightarrow\{1,\dots,n\}$ 
induce inclusions of groups $F_1\to F_n$
which express $F_n$ as a coproduct in $\T^{\rm op}_{\gr}$ of $n$ copies of
$F_1$; it follows that in the opposite category $\T_{\gr}$ the object
$F_n$ is the product of $n$ copies of $F_1$. 
Suppose that $ G$ is an 
arbitrary group. We can define a functor
$$ A_G\colon \T_\gr\lra\Spaces,\ \ \ \  
F_n\longmapsto \Hom_\gr(F_n,G).$$

It is clear that $A_G$ is product-preserving, 
and so $A_G$ is a strict $\T_\gr$-algebra. 
One can check that the converse is also true: 
any strict $\T_\gr$-algebra $A$ defines a group 
structure on the space $A(F_1)$. This is not surprising,
since (by Yoneda's lemma) the maps 
$F_n\to F_1$ in $\T_\gr$ correspond exactly 
to all of the ways of taking $n$ elements in 
a group and combining them with the available
operations to obtain a single element of the group. 
The composition in $\T_\gr$ gives identities 
between composites of these multivariable operations. 
A set which possesses such operations satisfying the 
appropriate identities is exactly a group.

Lawvere \cite{lavwere} showed that strict algebras 
can be used in this way to describe a wide 
class of algebraic structures, including,
besides groups, monoids, nilpotent and solvable 
groups of any fixed class, rings, Lie algebras etc. 
As  the example
above suggests,
 the existence of free objects is 
essential in order to get such a description. 

The language of algebraic theories proved to be 
equally convenient for describing various homotopy
invariant structures on spaces. However, in order
to allow for homotopy input one needs to relax 
the definition of a strict algebra.
Suppose that $\T$ is an 
algebraic theory with objects $T_n$,
$n\ge 0$. The expression of $T_n$ as a product of $n$ copies of $T_1$
gives projection maps
$$
     p^n_k:T_n\lra T_1,\ \ \ \    1\le k\le n.
$$

\numbereddemo{Definition}
\label{HOMOTOPYALG}
Suppose that $\T$ is an algebraic theory. 
A functor\break$X:\T\ra \Spaces$ is said to be a {\it homotopy}
$\T$-{\it algebra} if $X$ preserves products 
up to weak equivalence, i.e., if $X(T_0)$ is weakly contractible 
and for each $n\ge 1$ the product map
$$ \prod_{k=1}^{n}X(p^n_k)\colon X(T_n)
\mapright{20}{ } X(T_1)^n $$  
is a weak equivalence.
\enddemo

A homotopy $\T$-algebra structure on 
a space $Y$ is a homotopy
$\T$-algebra $X$ together with a weak
equivalence $X(1)\simeq Y$.  We can now state 
our main result which  implies that it is always
possible to pass from a homotopy $\T$-algebra 
structure on a space $Y$ to a strict $\T$-algebra 
structure on a space weakly equivalent to $Y$.

\proclaim{Theorem}
\label{MAIN}
Let $\T$ be an algebraic theory{\rm .} For any homotopy\break
$ \T$\/{\rm -}\/algebra  $X$ there exists a weak equivalence 
$ X\simeq LX$ such that $ LX$ is a strict
$ \T$\/{\rm -}\/algebra{\rm .}
\endproclaim

We will actually prove a somewhat stronger statement 
(\ref{NEWMAIN})
expressing the relationship of homotopy and strict 
$ \T$-algebras as a Quillen equivalence of  
model categories. In particular, the weak equivalence
in the theorem above respects the homotopy $\T$-algebra 
structures on both objects involved.

Theorem \ref{MAIN} gives a rigidifying result for 
homotopy algebras, but the following corollary 
shows that it is also of consequence for strict 
algebras.

\proclaim{{C}orollary}\label{FAY}
Let $ F\colon\Spaces\ra\Spaces$ be a functor
preserving weak equivalences and preserving products
up to weak equivalence{\rm .}
If $ Y$ is a space with a strict $ \T$\/{\rm -}\/algebra 
structure for some algebraic 
theory $ \T$ then $ F(Y)$ is weakly 
equivalent to a space with a strict 
$ \T$\/{\rm -}\/algebra structure{\rm .}
\endproclaim
 
Indeed, the assumptions on the functor $ F$
imply that for any strict\break $ \T$-algebra $ A$
the composition $F\!\circ\! A\colon\T\ra\Spaces$ 
is a homotopy $ \T$-algebra. Therefore, the statement
follows immediately from Theorem~\ref{MAIN}. 

Examples of functors for which Corollary~\ref{FAY}
holds include localization functors \cite{dror}
and Bousfield-Kan completion functors 
\cite{bousfield-kan}.

\numbereddemo{Note}
Although we define  an algebraic theory $\T$
as a discrete category (\ref{STRICTALG}) all 
statements of this paper remain valid 
also if we assume that $ \T$ is a simplicial 
category (and thus strict, and homotopy 
$ \T$-algebras are simplicial functors). 
The proofs in this case require at most 
minor changes. 
\enddemo

{\it Relationship to previous results}.
The notion of a homotopy algebra is inspired by 
$\Gamma $-spaces of Segal \cite{segal}. We note 
however that the indexing category $ \Gamma^{\rm op}$
which Segal uses is not an algebraic theory, but 
falls into a more general class of semi-theories:

\numbereddemo{Definition}
A {\it semi\/{\rm -}\/theory} $ \bo C$ is a small category with objects 
$C_0, C_1, \dots$ and such that for every $ n\ge 1$ 
there is a fixed set of morphisms $$ p^n_k\in 
\Hom_{\bo C}(C_n, C_1), \qquad
  1\leq k\leq n.$$
A functor $ X\colon {\bo C}\ra \Spaces$ is a {\it homotopy} 
(resp.\ {\it strict}\/)
$ \bo C$-{\it algebra} if $X(C_0)$ is weakly contractible 
(resp.\ $X(C_0)=\ast$)
and for $ n\ge 1$ the product  map
$$ \prod_{k=1}^{n}X(p^n_k)\colon X(C_n)\ra X(C_1)^n$$
is a weak equivalence 
(resp.\ an isomorphism).
\enddemo 

Segal proved that giving a homotopy 
$ \Gamma^{\rm op}$-algebra 
$ X$ amounts to providing the space $ X(C_1)$
with a structure of an infinite loop space up to group 
completion. Other examples of applications 
of homotopy algebras over semi-theories include 
characterization  of $n$-fold loop spaces 
\cite{bousfield} and generalized 
Eilenberg-Mac Lane spaces \cite{badzioch}.

In the study of strict algebras the passage 
from algebraic theories to semi-theories brings 
nothing new. In fact, for any semi-theory $ \C$
one can find an algebraic theory $ \bar\C$ such that 
the categories of strict $ \C$- and $ \bar\C$-algebras
are isomorphic. A result of this kind holds 
also for homotopy algebras, but in that case 
the construction of an algebraic theory associated 
to a semi-theory is more complicated \cite{badziochII}. 
It can be shown however, that many interesting homotopy
structures on spaces can be described directly as homotopy 
algebras over an algebraic theory. For example,
loop spaces can be characterized as homotopy algebras
over the theory $ \T_\gr$  which we have already mentioned. 
Also,\break $ A_\infty$-spaces and generalized 
Eilenberg-Mac Lane spaces can be viewed as homotopy
algebras over algebraic theories $ \T_{\bo{M}}$
and $ \T_{\bo{AbM}}$ such that the corresponding
strict algebras describe respectively monoids and 
abelian monoids. In each of these cases 
Theorem~\ref{MAIN} recovers the results mentioned
at the beginning of this paper.

\vglue12pt {\it Organization of the paper}.
In Section 2 we state some basic properties of algebraic 
theories and strict algebras.
Then, in Section 3, we recall two standard 
model category structures defined on 
a category of functors with values in $ \Spaces$. 
Section 4 contains some remarks on function 
complexes in model categories. 
In Section 5 we put a model category structure
on the category of strict algebras. Also, we  
describe a model category expressing  
the homotopy theory of homotopy $ \T$-algebras. 
Finally, in Section~6,  Theorem \ref{MAIN} is 
restated in the language of model categories and
proved in that form. 
 
\numbereddemo{Notation}
\label{NOTATION}
(i)\ This paper is written simplicially: by $ \Spaces$
we denote the category of simplicial sets. 
Consequently, by  `space' we always mean a simplicial set.
\vglue4pt
(ii)\  We use extensively the language of
model categories of Quillen. Our main 
references for this subject are
\cite{hirschhorn} and \cite{goerss-jardine}.
\vglue4pt
(iii)\  Given a category $ \bo M$ we will denote by $ s\bo M$
 the category of simplicial objects 
in $ \bo M$, that is, the category of functors 
$ {\bo \Delta}^{\rm op}\ra \bo M$. If $ \bo M$ is a model 
category then by the model category structure on $ s\bo M$
we will always understand Reedy model category structure
\cite[Ch.\ 16]{hirschhorn}, where weak equivalences are
objectwise weak equivalences while fibrations and 
cofibrations are defined using matching and latching 
objects. 
\vglue4pt
(iv)\  If $ \bo M$ is a simplicial model category 
then so is $ s\bo M$. In this case we have the geometric
realization functor 
$$ |-|\colon s{\bo M}\lra \bo M$$
defined by the coequalizer diagram: \pagebreak
$$ 
\coprod_{\phi\colon{\bo n}\ra {\bo m}}
X_m\otimes\Delta[n]\hskip4pt\raise2pt\hbox{$\lrar$}\hskip-23pt\lower2pt\hbox{$\lrar$} \coprod_{\bo
n}X_n\otimes\Delta[n]
 \lrar |X_\bullet| 
 $$
for $ X_\bullet\in s\bo M$ 
(see \cite[VII.3]{goerss-jardine}).
If $ \bo M = \Spaces$ and so $ s\bo M$ is the category 
of bisimplicial sets then $ |X_\bullet|$ is just 
the diagonal of $ X_\bullet$.
\enddemo
 
{\it Acknowledgement}.
This paper is a version of my Ph.D. thesis completed 
at the University of Notre Dame. I want to express my 
gratitude to my thesis advisor W. G. Dwyer who generously 
shared with me his ideas and knowledge during the time 
I was working on this project. In particular, he brought 
to my attention localization of model categories which 
proved to be the crucial tool used here. Also, his comments 
helped to improve the presentation of this paper.

\section{Algebraic theories}

We start with a brief review of algebraic
theories and their strict algebras.  
For a detailed exposition   
we refer to \cite[\S 3]{borceuxII}.
 
Let  $ \T$ be an algebraic theory and 
$\alg$ the category of strict $ \T$-algebras.
We have the forgetful functor 
$$ U_{\T}\colon{\alg}\ra{\Spaces},\ \ \ \  U_{\T} (A):=A(T_1).$$
 It is in fact a half of an adjoint pair:

\proclaimtitle{\cite[2.3]{schwede}}
\proclaim{Proposition}
\label{prop23}
The functor $ U_{\T}$ has a left adjoint 
$F_{\T}\colon\Spaces\break\ra\alg${\rm .} If $ Y\in\Spaces$
then
$$F_\T(Y)(T_1)=\coprod_{n\geq 0}\Hom_\T(T_n,T_1)
\times Y^n/\sim$$
where the identifications come from the set operations
present in any algebraic theory{\rm .} 
\endproclaim

The functor $ F_{\T}$ 
will be called the free $ \T$-algebra functor.

Let $ \st$ be the category of all simplicial functors 
$ \T\ra\Spaces$.
We will often identify the category $ \alg$  
with a full subcategory of $ \st$. Using 
this identification we get:

\proclaim{Proposition}
\label{ALGLIMITS}
The category $ \alg$ is complete and the limits are
computed objectwise{\rm . } 
\endproclaim

\demo{Proof}
All limits in $ \st$ exist and are computed objectwise, so 
it is enough to notice that a limit of product-preserving functors also preserves products.
\enddemo 

Let $ J_{\T}\colon\alg\ra\st$ denote 
the embedding of categories.
Proposition ~\ref{ALGLIMITS} immediately implies:

\proclaim{{C}orollary}\label{PRESLIMITS}
The functor $ J_{\T}$ preserves 
limits{\rm .}
\endproclaim

The following fact shows that $ \alg$
is a reflective subcategory of $ \st$.

\proclaim{Proposition}
\label{LEFTADJTOINCLUSION}
There exists a functor $ K_{\T}\colon\st\ra\alg$
left adjoint to $ J_{\T}${\rm .}
\endproclaim

\demo{Proof}
We use the adjoint functor theorem 
\cite[Thm.\ 2, p.\ 117]{maclane}.
By (\ref{ALGLIMITS}) and (\ref{PRESLIMITS}) it is 
enough to check the solution set condition. 

Let $ f\colon X\ra A$ be a morphism in $ \st$
such that $ A\in\alg$. For $ n\geq 0$ let
$$f_n\colon X({T}_n)\ra A({T}_n) $$
denote the restriction of $ f$ to $ T_n$.
By (\ref{prop23}) the map 
$ f_1$ has a left adjoint
$$g\colon F_\T(X(T_1))\ra A.$$ 
Define $ M_f$  to be the image of $ g$: 
$$M_f(T_n):=
{\rm im}\left(F_\T(X(T_1))(T_n)\mapright{15}{g_n} 
A(T_n\right)).$$
Since $g$ is a map of strict $ \T$-algebras  
we have $ g_n = (g_1)^n$
and it follows that  $M_f$ is also a strict algebra.
We claim that there exists a morphism 
$\bar{f}\colon X\ra M_f$ such that 
the following diagram commutes:
\figin{p6}{1000}%
Indeed, it is enough to show that  for $ n\geq 1$
the image of $ f_n$ is contained in 
$ M_f(T_n)$. 
For $ n=1$ this follows directly from
the definition of $ M_f$. For $ n>1$,
since both $ A$ and $ M_f$
are strict $ \T$-algebras,  we have
$${\rm im}(f_n)\subseteq {\rm im}(f_1)^n$$  
and
$$ {\rm im}(f_1)^n\subseteq (M_f(T_1))^n
\cong M_f(T_n).$$
Therefore ${\rm im}(f_n)\subseteq M_f(T_n)$
as claimed.  

Let $ \lambda$ denote the cardinality of 
the set of simplices of the space
$\coprod_{n}X(T_n)\times \Hom_\T(T_n,T_1)$.
Since $F_\T(X(T_1))$ maps onto $M_f$, from the description of
the functor $ F_\T$ (\ref{prop23}) 
we get that the cardinality of the set of 
simplices of $M_f(T_1)$ cannot exceed  
$ \lambda$. Hence the solution set for $ X$
can be chosen to consist of 
the representatives of isomorphism classes of 
these strict $\T$-algebras $ B$ for which the set of
simplices of $ B(T_1)$ has cardinality not 
greater than $ \lambda$.
\enddemo

The category $ \st$ is cocomplete;  therefore
(\ref{LEFTADJTOINCLUSION}) implies:

\proclaim{{C}orollary} The category $ \alg$ is cocomplete{\rm .}
\endproclaim

Colimits in $ \st$ are computed objectwise
but this is not true in general in $ \alg$. 
However one has the following:

\proclaimtitle{\cite[3.4.2]{borceuxII}}
\proclaim{Proposition}
\label{FILTEREDCOLIMS}
Filtered colimits in $ \alg$ are computed 
objectwise{\rm .} In particular the inclusion functor
$ J_\T\colon \alg\ra\st$ preserves filtered 
colimits{\rm .}
\endproclaim

In Section \ref{STRICT AND HOMOTOPYALGEBRAS}
we will need the following observation which is 
an easy consequence of the adjointness of $F_\T$ and $U_\T$.

\proclaim{Lemma}
\label{K_T OF COPRODUCTS}
For $  1\leq i\leq m$
let $[n_i]=\{1,\dots,n_i\}$ be a discrete 
simplicial set and let
$$ \kappa\colon \coprod_{i=1}^m F_\T ([n_i])\ra
F_\T\left(\coprod_{i=1}^m [n_i]\right)$$
be the map in $\st$ induced 
by inclusions $[n_i]\ra \coprod_{i=1}^m [n_i]${\rm .}
Then $ K_\T(\kappa)$ is an isomorphism in 
$ \alg${\rm .} 
\endproclaim

Observe that if $ A$ is a strict algebra then 
$ K_\T A\cong A$. Therefore we get:

\proclaim{{C}orollary}\label{COR TO K_T OF COPRODUCTS}
If $ [n_i]$ are simplicial sets as above{\rm ,} 
$ 1\leq i \leq m$, then  
$$ K_\T\left(\coprod_{i=1}^m F_\T ([n_i])\right)\cong
F_\T\left(\coprod_{i=1}^m [n_i]\right).$$
\endproclaim  

\numbereddemo{{R}emark}
\label{F_N}
 For $ [n]=\{1,\dots,n\}$ we will denote the functor
$ F_\T([n])$ by $F_n$. Notice that $ F_n$ can be 
described as a functor corepresented by $ T_n \in \T$:
$$ F_n(T_m)=\Hom_{\T}(T_n, T_m). $$
Indeed, for every strict $ \T$-algebra $ A\in\alg$ we have 
$$\Hom_{\alg}(\Hom_{\T}(T_n,-), A)\cong A(T_n)\cong
A(T_1)^n \cong \Hom_{\Spaces}([n], A(T_1))$$
where the first isomorphism comes from
Yoneda's lemma \cite[p.61]{maclane}. But $ F_\T$ is left 
adjoint to the forgetful functor $ U_\T$ 
(\ref{prop23}), and so we have an isomorphism 
$ \Hom_{\Spaces}([n], A(T_1))\cong \Hom_\alg(F_n, A)$. 
It follows that $ \Hom_{\T}(T_n, -)$ must be isomorphic 
to $ F_n$.
\enddemo

\section{Model category structures 
on categories of diagrams}

We recall here two standard model category 
structures defined on a category of diagrams 
of spaces. The model categories describing 
the homotopy theories of strict and homotopy 
$ \T$-algebras (\S \ref{MODEL CAT FOR ALGEBRAS})
will be derived  from $\stf$.
The properties of $ \stc$ on the other hand 
will allow us to avoid the trouble 
of working with homotopy function complexes 
(\S \ref{FUN COMPLEX}) as shown in 
(\ref{S-EQUIV. VIA SIMPL. COMPLEXES }).

Let $ \C$ be a small  category and 
let $ \sC$ denote the category of all functors 
$\C\ra\Spaces$.

\numbereddemo{Notation}
\label{MOD CAT OF DIAGRAMS}
Let $ \scf$ and $ \scc$ denote the category $\sC$
together with a choice of three classes of morphisms:
\begin{itemize}
\item{$\scf$}
\begin{itemize}
\item{weak equivalences := objectwise weak equivalences}
\item{fibrations := objectwise fibrations}
\item{cofibrations := morphisms with the left lifting property
with respect to all fibrations which are weak equivalences}
\end{itemize}
\item{$\scc$}
\begin{itemize}
\item{weak equivalences := objectwise weak equivalences}
\item{cofibrations := objectwise cofibrations}
\item{fibrations := morphisms with the right lifting property
with respect to all cofibrations which are weak equivalences}
\end{itemize}
\end{itemize}

\enddemo

\proclaimtitle{\cite[IX 1.4, VIII 2.4]{goerss-jardine}}
\proclaim{Theorem}
\label{MODEL STRUCTURES}
Both $ \scf$ and $\scc$ are simplicial
model categories{\rm .} In each case the simplicial structure is given by 
$$ (X\otimes K)(c)= X(c)\times K$$
for any $X\in\sC ${\rm ,} $ c\in\C$ and a simplicial set $ K${\rm .}
\endproclaim

Directly from (\ref{MOD CAT OF DIAGRAMS}) 
one gets that every object 
of $ \scc$ is cofibrant. In the remainder of this section 
we describe a canonical construction of a 
cofibrant replacement of a diagram of spaces with 
respect to $ \scf$ model category structure.

For a category $\C$ as above let $\Cdi$ denote 
the category with the same objects as $ \C$ and with no 
nonidentity morphisms. 
The following is readily verified.

\proclaim{Proposition}
\label{ADJ TO DISCRETE}
The forgetful functor $ U\colon\sC\ra\sCdi$ has a left 
adjoint $ F\colon\sCdi\ra\sC$ given by 
$$F(X):= \coprod_{c\in\C}F_c\otimes X(c)$$
where $ F_c\in \sC$ is a functor such that 
$F_c(d):=\Hom_{\C}(c,d)${\rm .} 
\endproclaim
 
Let $\eta \colon Y\ra UFY$ and 
$\varepsilon\colon FUX\ra X$ denote the unit 
and the counit of this adjunction. 
As for any pair of adjoint functors
the composition $$FU\colon\sC\ra\sC$$ defines a cotriple
(comonad) \cite[p.135]{maclane} with the structure maps 
$$ \varepsilon\colon FUX\ra X\ \ \ {\rm and}\ \ \ 
\delta\colon FUX\ra (FU)^2X$$
where $ \delta=F\eta U$.

\numbereddemo{Definition}
\label{RESOLUTION}
For $ X\in\sC$ the {\it standard simplicial resolution} of 
$X$ is a simplicial object $ FU_\bullet X\in s\sC$
which in the dimension $k$ consists of a diagram 
$ FU_k X:=(FU)^{k+1} X$. Face and 
degeneracy operators of $ FU_\bullet X$ are given by 
$$(FU_k X\mapright{15}{d_i}FU_{k-1}X):=
\left((FU)^{k+1}X 
\mapright{25}{(FU)^{i}\varepsilon (FU)^{k-i}}
(FU)^k X\right)$$
and
$$
(FU_k X\mapright{15}{s_i}FU_{k+1}X):=
\left((FU)^{k+1}X 
\mapright{25}{(FU)^{i}\delta (FU)^{k-i}} 
(FU)^{k+2} X\right).
$$
\enddemo 

If we regard $ X$ as a constant simplicial object 
we can define a simplicial map 
$ \varphi\colon FU_\bullet X\ra X$
$$(FU_k X\mapright{15}{\varphi_k} X):= 
((FU)^{k+1}X\mapright{20}{\ \varepsilon^{k+1}} X).$$ 
Let $ |\varphi|\colon |FU_\bullet X|\ra X$ be the 
geometric realization of $ \varphi$ 
(\ref{NOTATION}) taken with respect to the simplicial 
structure as in (\ref{MODEL STRUCTURES}).

\proclaim{Proposition}
\label{AUGMENTATION IS W.E.}
The map $|\varphi|$ is a weak equivalence{\rm .}
\endproclaim

\demo{Proof}
If  $Y_\bullet$ is a simplicial object in $\sC$ then
its realization can be computed objectwise:
$|Y_\bullet|(c)=|Y_\bullet (c)|$ where the space on the 
left-hand side is the realization (diagonal) of the
bisimplicial set $Y_\bullet (c)$. Therefore, it suffices 
to show that 
$|\varphi^c|\colon |FU_\bullet X(c)|\ra|X(c)|$
is a weak equivalence of spaces for all $ c\in \C$.
One can check that the realization of the map 
$\phi^c\colon X(c)\ra FUX(c)$ given by 
$$(X(c)\mapright{15}{\phi^c_k} FU_k X(c)):= 
(X(c)\mapright{15}{\eta^{k+1}}(FU)^{k+1} X(c))$$
is a homotopy inverse for $ |\varphi^c|$.  
\enddemo

We claim that $ |FU_\bullet X|$ is a cofibrant 
replacement for $X$. In view of 
(\ref{AUGMENTATION IS W.E.}) it remains to show 
that $ |FU_\bullet X|$ is a cofibrant object of
$ \scf$.  

\proclaim{Proposition}
\label{REEDY COFIBRANCY}
For every $ X\in \sC $ the resolution $FU_\bullet X$
is a Reedy cofibrant object in $s\scf${\rm .}
\endproclaim

\demo{Proof}
Let $ {\bo\Delta}_{+}^{\rm op}$ denote the subcategory
of $ {\bo \Delta}^{\rm op}$ generated by all degeneracy 
maps $ s_i$ and positive face maps $ d_j$, $ j>0$.
The category $ {\bo\Delta}_{+}^{\rm op}$ is a Reedy 
category \cite[16.1.2]{hirschhorn} with the direct 
subcategory of $ {\bo\Delta}_{+}^{\rm op}$
generated by degeneracy maps and the inverse 
subcategory generated by positive face 
maps.  For a model category $ \bo M$
define $ s_{\!\!+}\bo M$ to be the category of functors 
${\bo\Delta}_{+}^{\rm op}\ra \bo M$ with the Reedy model 
category structure. 

The embedding of categories
$ {\bo\Delta}_{+}^{\rm op}\hookrightarrow
{\bo\Delta}^{\rm op}$ defines a functor
$ {\rm\bo r}\colon s{\bo M}\ra s_{\!\!+}\bo M$. 
Since all degeneracy maps of 
$ {\bo\Delta}^{\rm op}$ are contained 
in $ {\bo\Delta}_{+}^{\rm op}$
we get that $ X_\bullet\in s{\bo M}$ is 
cofibrant if and only if $ {\rm\bo r}(X)$ is cofibrant in 
$ s_{\!\!+}\bo M$. Take $ {\bo M} = \scf$.
To prove the statement of 
the proposition it is then enough to
show that $ {\rm\bo r}(FU_\bullet X)$ is a cofibrant object 
in $ s_{\!\!+}\scf$. 
\vglue4pt
Since the category $ \sCdi$  is isomorphic
to the product $ \prod_{c\in \C}\Spaces$ it is 
a model category with weak equivalences, 
fibrations and cofibrations defined to be objectwise 
weak equivalences, fibrations and cofibrations.
The adjoint pair of functors $ (F,U)$ becomes 
then  a Quillen pair between $ \sCdi$ and $ \scf$.
It follows that the induced functors 
$$
F\colon s_{\!\!+}\sCdi\hskip4pt \raise2pt\hbox{$\lrar$} \hskip-24pt 
\lower2pt\hbox{$\llar$}\hskip4pt s_{\!\!+}\scf\colon U  
 $$
also form a Quillen pair with respect 
to Reedy model category 
structures\break \cite[16.11.1]{hirschhorn}. 
From Definition \ref{RESOLUTION} we see 
that the object ${\rm\bo r}(FU_\bullet X)$ 
is in the image of the functor $ F$, hence 
it suffices to show that every 
object of $s_{\!\!+}\sCdi$ is cofibrant.   
This is however an immediate consequence of 
the fact that every object of $s_{\!\!+}\Spaces$ is 
Reedy cofibrant. The proof of this last 
statement is the same as the proof that every 
object is cofibrant in the Reedy model 
category structure  on $s\Spaces$ -- the category
of bisimplicial sets (see \cite[16.7.8]{hirschhorn}).  
\enddemo

Geometric realization of a Reedy cofibrant object 
is cofibrant 
\cite[VII 3.6]{goerss-jardine},
hence (\ref{REEDY COFIBRANCY}) implies

\proclaim{{C}orollary}\label{COFIBRANCY}
The diagram $ |FU_\bullet X|$ is a 
cofibrant object of $\scf${\rm .}
\endproclaim

\section{Function complexes}
\label{FUN COMPLEX}

In Section \ref{MODEL CAT FOR ALGEBRAS} we will 
introduce a model category for homotopy 
$ \T$-algebras. In preparation for that  
we recall here some properties of function complexes 
in model categories.

Let $ \M$ be a simplicial model category and
let  $ \HoM$ denote its homotopy category. In 
\cite{dwyer-kanII} Dwyer and Kan showed that
for any $ X, Y\in \M $ the set of morphisms 
$\Hom_\HoM (X,Y)$ can be replaced by a richer 
structure of a homotopy
function complex, that is a simplicial set 
$ \Map_\M(X,Y)$ such that $ \pi_0\Map_M (X, Y)\cong \Hom_\HoM (X,Y)$. 
Moreover the following holds:

\begin{itemize} 
\item[(i)] $ \Map_M(X,Y)$ preserves weak 
equivalences:
if $X\simeq X'$ and $Y\simeq Y'$ then 
$ \Map_\M(X,Y)\simeq \Map_\M(X',Y')$.
\item[(ii)] The homotopy type of $ \Map_\M(X,Y)$
depends only on the class of weak equivalences of 
$ \M$: if $\M'$ is a model category with the same 
underlying category as $\M$ and with the same
class of weak equivalences then
$ \Map_\M(X,Y)\cong\Map_{\M'}(X,Y)$ for 
all $ X, Y$ in $\M$. 
\item[(iii)] A morphism $ f\colon X\ra X'$ is a weak 
equivalence if and only if the induced map 
$ f_\ast\colon\Map_\M(X',Y)\ra\Map_\M(X,Y)$ 
is a weak equivalence of simplicial sets for all
$ Y\in\M$.
\item[(iv)] If $ \M$ is a simplicial model category,
$ X$ is cofibrant and  $ Y$ is fibrant
then $\Map_\M(X,Y)$ is weakly equivalent 
to the simplicial set $\map_\M(X,Y)$ where 
$\map_\M(X,Y)_k = \Hom_\M(X\otimes\Delta[k], Y)$.
We will call $ \map_\M(X,Y)$ the simplicial 
function complex of $ X$ and $ Y$.   
\end{itemize}

\numbereddemo{Note}
As a consequence of (ii) we do not need to 
distinguish between homotopy function complexes 
taken with respect to $ \scf$ and $ \scc$ model 
category structures. Hence from now on  
$ \Map(-,-)$ will stand for a homotopy function 
complex in any of these model categories. 
Similarly, since simplicial function complexes
(see (iv)) are defined using only the  simplicial structure
of a model category, they are the same in $\scf$ and 
$\scc$. We will denote them by $ \map(-,-)$. 
\enddemo

In Section \ref{STRICT AND HOMOTOPYALGEBRAS}
we will refer to the following property of 
simplicial function complexes.

\proclaim{Proposition}
\label{W.E. FOR SIMPLICIAL OBJECTS}
Let $ \M$ be a simplicial model category and let 
$ Y\in \M$ be a fibrant object{\rm .} Assume that 
$ f\colon X_\bullet\ra X_\bullet'$ is a map 
of Reedy cofibrant objects in $ s\M$ such that 
$$f_n^\ast \colon \map_\M({X_n}', Y)
\,\mapright{15}{ }\,
\map_\M(X_n, Y) $$
is a weak equivalence for all $ n\geq 0${\rm .} 
Then the geometric realization of the map $ f$ induces 
a weak equivalence
$$|f|^\ast\colon \map_\M(|X_\bullet'|,Y)
\,\mapright{15}{\simeq}\, \map_\M(|X_\bullet|, Y).$$
\endproclaim

\demo{Proof}
We have a commutative diagram
\figin{p11}{800}
The map $ \Phi_X^\ast$ is induced by the Bousfield-Kan
map $ \Phi_X\colon {\rm hocolim\,}_{\Delta^{\rm op}}X_\bullet
\ra |X_\bullet|$ \cite[19.6.3]{hirschhorn} 
(and similarly for $ \Phi_{X'}^\ast$), 
while  $ \Psi_X$ and $\Psi_{X'}$ 
are the isomorphisms of simplicial 
sets described in \cite[19.1.12]{hirschhorn}.
Since $X_\bullet$ and $X_\bullet '$ are Reedy 
cofibrant it follows
from \cite[19.6.4]{hirschhorn} that 
$ \Phi_X$ and $\Phi_{X'}$ are weak equivalences 
and so are the maps they induce on simplicial 
function complexes.
Fibrancy of $ Y$ implies on the other hand, that $ \map(-,Y)$
is always a Kan complex. Therefore, by our assumption 
on $ f$ the bottom map $f^\ast$ is a weak equivalence
(see \cite[19.4.3]{hirschhorn}),
and hence so is the top map. 
\enddemo

\numbereddemo{Note}
\label{4.2 FOR SCC}
Suppose than in Proposition 
\ref{W.E. FOR SIMPLICIAL OBJECTS} we have 
$ \M=\scc$ (\ref{MOD CAT OF DIAGRAMS}). Then 
the assumption that $ X_\bullet$ and 
$ X_\bullet'$ are Reedy cofibrant is always satisfied. 
Indeed, since cofibrations in $ \scc$ are defined 
objectwise, $ X_\bullet\in s\scc$ is Reedy cofibrant
if and only if $ X_\bullet(c)$ is a cofibrant bisimplicial set 
for all $ c\in\C$. This last condition however 
always holds since all bisimplicial sets are cofibrant
in the Reedy model category structure on $ s\Spaces$
\cite[16.7.8]{hirschhorn}.
\enddemo
 
\section{Model category for homotopy $ \T$-algebras}
\label{MODEL CAT FOR ALGEBRAS}

Let $ \T$ be an algebraic theory.
Recall that by (\ref{LEFTADJTOINCLUSION})
there is an adjoint pair of functors 
$$
K_\T\colon \st\hskip3pt \raise2pt\hbox{$\lrar$}\hskip-24pt\lower2pt\hbox{$\llar$}\hskip3pt \alg\colon J_\T  $$
where $ J_\T$ is the inclusion onto a subcategory.  
This adjunction can be used to put a model category
structure on $ \alg$:  

 \proclaimtitle{\cite[3.1]{schwede}}
\proclaim{Theorem}
\label{MODEL STR ON ALG}
The category $ \alg$ is a model category 
with  weak equivalences and fibrations
defined  as  objectwise  weak equivalences and fibrations.
Then the adjoint pair $(K_\T, J_\T)$ becomes a Quillen 
pair between $ \stf$ and $ \alg$. 
\endproclaim

Our main goal in this 
section is to construct a model category $ \lst$
which reflects the homotopy theory of homotopy 
$ \T$-algebras.

Recall (\ref{F_N}) that for each 
$ n\geq 0$ there is a functor
$ F_n\in \st$ given by $ F_n(T_m):= \Hom_{\T}(T_n, T_m)$
(in fact, since $ T_m\cong T_1^m$ we get that 
$ F_n\in \alg$). For $ n\geq 1$ the projections
$ p^n_k$ induce maps $ p_n\colon \coprod_n F_1\ra F_n$.
We define also $ p_0\colon \coprod_0 F_1\ra F_0$
to be the unique map from the diagram of 
empty spaces to $F_0$
(alternatively, the morphism $ p_n$ can be
described as the map induced by inclusions of sets 
$ [1]\hookrightarrow [n]$ as in Lemma 
\ref{K_T OF COPRODUCTS}).
Let $ S:=\{p_0, p_1,\dots \}$.

\numbereddemo{Definition}
\label{LOC EQ}
An {object} $ Z\in \stf$ is $ S$-\/{\it local} if it is 
fibrant and if for each $n\geq 0$ the map of
homotopy function complexes 
$$p_n^\ast\colon\Map(F_n,Z)\ra\Map(\coprod_n{F_1},Z)$$
is a weak equivalence.

A {morphism} $ f\colon X\ra Y$ in $ \stf$ is an
$S$-{\it local equivalence} if the induced map 
$$f^\ast\colon\Map(Y,Z)\ra\Map(X,Z) $$
is a weak equivalence for every $ S$-local object $ Z$.  
\enddemo

\numbereddemo{Note} 
\label{USE SIMPL COMPLEX}
Since both $\coprod_n{F_1}$ and $F_n$ are cofibrant 
and $Z$ is fibrant in $\stf$ the map $p_n^\ast$ in the definition 
of $S$-local objects above can be replaced by the map 
of simplicial function complexes (\S \ref{FUN COMPLEX}(iv)):
$$p_n^\ast\colon\map(F_n,Z)\ra\map\left(\coprod_n{F_1},Z\right).$$ 
\enddemo

\proclaim{Proposition}
\label{LOCAL STRUCTURE}
Let $ \lst$ denote the category $ \st$ with three
distinguished classes of morphisms\/{\rm :}
\begin{itemize}
\item[--]weak equivalences $:=  S$-local equivalences
\item[--]cofibrations $:=$ cofibrations in $ \stf$
\item[--]fibrations $:=$ maps with the right lifting 
property with respect to all cofibrations which are
weak equivalences
\end{itemize} 
Then $ \lst$ is a simplicial model category with the same
simplicial structure as $ \stf$ 
{\rm (\ref{MODEL STRUCTURES}).}
\endproclaim

\demo{Proof}
This is a consequence of a general result 
\cite[4.1.1]{hirschhorn} which  proves the existence 
of left Bousfield localizations for a broad 
class of model categories. The category $ \stf$
satisfies the assumptions of that theorem by 
\cite[4.1.4]{hirschhorn} and \cite[4.1.5]{hirschhorn}
and the model category structure on $ \lst$
is obtained by localizing $ \stf$ with respect to 
the set $ S$. 

\enddemo

Within the model category $ \lst$, homotopy 
$\T$-algebras 
can be characterized as follows:

\proclaim{Proposition}
\label{LOCAL FIBRANTS}
An object $ Z\in\lst$ is fibrant if and only if
it is a homotopy ${\bf T}$\/{\rm -}\/algebra{\rm ,}
fibrant as an object of $ \stf$.
\endproclaim 

\demo{Proof}
By \cite[3.5.1]{hirschhorn}
fibrant objects of $ \lst$ are exactly the $ S$-local
objects. 
Therefore, for any fibrant $ Z\in\lst$ 
the maps $p_n^\ast$ as in (\ref{USE SIMPL COMPLEX}) are 
weak equivalences. But for every $ n\geq 0$ we have
$ \map(\coprod_n{F_1},Z)\cong\prod_n{Z(T_1)}$
and $ \map(F_n,Z)\cong Z(T_n)$. It follows
that  $ Z$ is a homotopy $\T$-algebra. 
The proof of the other implication is similar.     
\enddemo

The next proposition is a consequence of 
(\ref{LOCAL FIBRANTS})
and \cite[3.3.12]{hirschhorn}.

\proclaim{Proposition}
\label{LOCAL W.E.}
If $ Z,Z'\in\st$ are homotopy $ \T$-algebras 
and $ f\colon Z\ra Z'$ is an $ S$-local weak 
equivalence then $ f$ is a weak equivalence in
$ \stf$ (i.e. an objectwise weak equivalence).
\endproclaim

Proposition \ref{LOCAL FIBRANTS} 
says that the homotopy category of 
$ \lst$ 
is equivalent to the category of homotopy  
$\T$-algebras with inverted $ S$-local equivalences.
From  (\ref{LOCAL W.E.}) we see however  that it 
amounts to inverting objectwise weak equivalences.
As a consequence we get:  

\proclaim{{C}orollary}\label{MODEL FOR HOMOTOPYALGEBRAS}
The homotopy theory of homotopy $\T$\/{\rm -}\/algebras
{\rm (}\/with objectwise weak equivalences\/{\rm )} is 
equivalent to the homotopy category of $ \lst${\rm .}
\endproclaim

The definition of $ S$-local equivalences we gave  
above (\ref{LOC EQ}) involves maps defined 
on homotopy function complexes 
$ \Map(-,Z)$. In practice it is more
convenient to work with simplicial function 
complexes $ \map(-,Z)$. Since we assume that $ Z$
is fibrant in $ \stf$ we get 
$ \Map(X,Z)\simeq\map(X,Z)$ whenever $ X\in\stf$
is a cofibrant object. However, the property that 
$ f\colon X\ra X'$ is an $ S$-local equivalence 
can be expressed in terms of simplicial function 
complexes even when $ X$ or $ X'$ is not cofibrant.

\proclaim{Proposition}
\label{S-EQUIV. VIA SIMPL. COMPLEXES }
A map $ f\colon X\ra X'$ is an $ S$\/{\rm -}\/local equivalence
if and only if for any homotopy $ \T$\/{\rm -}\/algebra $ \bar{Z}$ fibrant 
in $ \stc$ the induced map 
$$f^\ast\colon\map(X',\bar{Z})\ra\map(X,\bar{Z}) $$
is a weak equivalence of simplicial sets\/{\rm .}\/
\endproclaim

\demo{Proof}
Assume that $ f$ is an $ S$-local equivalence. 
Since fibrant objects of $ \stc$ are also fibrant in 
$ \stf$ the map $ f$ induces a weak equivalence on 
$ \Map(-,\bar{Z})$ for any $\bar{Z}$ as above. 
Moreover, all objects of $ \stc$ are cofibrant; hence 
we have $ \Map(Y,\bar{Z})\simeq \map(Y,\bar{Z})$ for any 
$ Y\in\stc$. It follows that the map
$ f^\ast\colon\map(X',\bar{Z})\ra\map(X',\bar{Z})$
must be a weak equivalence.

Conversely, assume that $f$ induces a weak equivalence
on $ \map(-,\bar{Z})$ for all homotopy $\T$-algebras 
$ \bar{Z}$,  fibrant in $ \stc$. Let $ Z$ be any 
homotopy $ \T$-algebra and let 
$ Z\mapright{15}{\simeq}\bar{Z}$ denote a fibrant 
replacement of $ Z$ in $ \stc$. Then 
$ \Map(-,Z)\simeq\Map(-,\bar{Z})\simeq\map(-,\bar{Z})$,
and so $f^\ast\colon\Map(X',Z)\ra\Map(X,Z)$ is a 
weak equivalence. 
Therefore $ f$ is an $ S$-local equivalence.
\enddemo

 \vglue-16pt

\section{Strict and homotopy $ \T$-algebras}
\label{STRICT AND HOMOTOPYALGEBRAS}
\vglue-8pt

We begin by proving some properties of the adjunction
$(J_\T,K_\T)$ (see (\ref{LEFTADJTOINCLUSION})).
Let $ A_\bullet$ be a simplicial object in $ \alg$.
Abusing   notation, by $|A_\bullet|$ we will always
denote the geometric realization of $ A_\bullet$
with respect to the simplicial structure in $ \stf$
(\ref{MODEL STRUCTURES}), that is, the geometric 
realization of $ J_\T A_\bullet\in\stf$.

\proclaim{Lemma}
\label{REALIZATION OF ALGEBRA}
If $ A_\bullet\in s\alg$ then $ |A_\bullet|\in\alg${\rm .}
\endproclaim

\demo{Proof}
We need to show that $ |A_\bullet|(T_0)\cong\ast$
and that for $ n>0 $ the projection maps
$ p_k^n\colon T_n\ra T_1$ induce isomorphisms 
$ |A_\bullet|(T_n)\cong(|A_\bullet|(T_1))^n$.
As we have already noted in the proof of 
(\ref{AUGMENTATION IS W.E.}) the realization
of $ A_\bullet$ can be computed objectwise:
$$|A_\bullet|(T_n)\cong |A_\bullet(T_n)|$$
where $|A_\bullet(T_n)|$ is the diagonal 
of the bisimplicial set $ A_\bullet(T_n)$.
Since $ A_m(T_0)=\ast$ for all $ m$,
we have $ |A_\bullet(T_0)|=\ast$. 
Moreover, since $ A_\bullet\in s\alg$, 
the projection maps induce
isomorphisms of bisimplicial sets 
$ A_\bullet(T_n)\mapright{15}{\cong}(A_\bullet(T_1))^n$
for\break $ n> 0$. 
Since the diagonal of a bisimplicial set commutes 
with products we get 
$ |A_\bullet(T_n)|\cong|A_\bullet(T_1)|^n$,
and it follows that $ |A_\bullet|$ is a strict
$ \T$-algebra. 
\enddemo

\proclaim{Lemma}
\label{REALIZATION OF K_T}
If $ X_\bullet$ is a simplicial object in $\stf$
then $ K_\T|X_\bullet|\cong|K_\T X_\bullet|${\rm .}
\endproclaim

\demo{Proof}
By (\ref{REALIZATION OF ALGEBRA}),
$|K_\T X_\bullet|$ is an object of $ \alg$.
Let $ \eta_\bullet\colon X_\bullet\ra K_\T X_\bullet$
($=J_\T K_\T X_\bullet$) be the unit of the 
adjunction $ (K_\T,J_\T)$. By properties of adjunction
there exists a map 
$ \theta\colon K_\T|X_\bullet|\ra |K_\T X_\bullet|$
such that the following diagram commutes:
\figin{p15}{1000}%
The map $ \eta_{|X_\bullet|}$ above 
is the unit of adjunction 
for $ |X_\bullet|$. It is enough to prove that 
$ |\eta_\bullet|^\ast\colon
\Hom(|K_\T X_\bullet|, A)\ra \Hom (|X_\bullet|, A)$ 
is an isomorphism for all 
$ A\in\alg$. Indeed, by properties of adjunction  
$ \eta_{|X_\bullet|}^\ast\colon \Hom(K_\T|X_\bullet|, A)
\ra\Hom(|X_\bullet|, A)$ is an isomorphism for any 
$ A$, so it will follow that
the map $ \theta^\ast\colon\Hom(K_\T|X_\bullet|,A)\ra
\Hom(|K_\T X_\bullet|, A)$ must be an isomorphism 
for all $ A\in \alg$, and hence $ \theta $ 
is an isomorphism. 
For a cosimplicial space $ Y_\bullet$ let
$ {\rm Tot}(Y_\bullet)$ denote the realization 
of $ Y_\bullet$ \cite[VIII.1 p.390]{goerss-jardine}.
We have an isomorphism
$$\Hom(|X_\bullet|,A)\cong
{\rm Tot}(\Hom(X_\bullet, A))_0 .$$ 
Therefore, it suffices to show that the realization of 
the map of cosimplicial spaces 
$ \eta_\bullet^\ast\colon\Hom(K_\T X_\bullet, A)
\ra \Hom(X_\bullet, A)$ is an isomorphism.
But this follows from the fact that
$ \eta_n\colon \Hom(K_\T X_n,A)\ra\Hom(X_n,A)$
is an isomorphism for any strict $ \T$-algebra $ A$
and for any $ n\geq 0$.  
\enddemo 

Next, recall that the adjunction $ (K_\T,J_\T)$ is 
a Quillen pair between\break $ \stf$ and $ \alg$
(\ref{MODEL STR ON ALG}).
The following fact shows that this property 
does not change if we replace the model category 
structure on $ \st$ by the $ S$-local structure.

\proclaim{Proposition}
The adjoint functors $(K_\T, J_\T) $ form a 
Quillen pair between the categories
$ \lst $ and $ \alg${\rm .}
\endproclaim

\demo{Proof}
Let $ p_n\colon\coprod_n F_1\ra F_n$ 
be the map given in Section~\ref{MODEL CAT FOR ALGEBRAS}. 
The model category $\lst$ is obtained by localizing 
$\stf$ with respect to all maps $p_n$. 
Thus the above statement follows from the observation  
that by (\ref{K_T OF COPRODUCTS}) $ K_\T(p_n)$  
is an isomorphism in $ \alg$ and 
from \cite[3.4.20]{hirschhorn}. 
\enddemo

Since $ \lst$ can serve as a model category 
for the homotopy theory of homotopy $ \T$-algebras
(\ref{MODEL FOR HOMOTOPYALGEBRAS}), out main 
Theorem \ref{MAIN} can be now restated as follows.

\proclaim{Theorem}
\label{NEWMAIN}
The Quillen pair of functors
$$ 
K_\T\colon \lst\hskip3pt\raise2pt\hbox{$\lrar$}\hskip-24pt\lower2pt\hbox{$\llar$}\hskip3pt \alg\colon J_\T
 $$
is a Quillen equivalence. 
\endproclaim

This is in turn a consequence of the following: 

\proclaim{Lemma}
\label{UNIT OF ADJUNCTION}
Let $ \eta_{X}\colon X\ra K_\T(X)\ (=J_\T K_\T(X))$
denote the unit of the adjunction $(K_\T,J_\T)${\rm .}
Then for every cofibrant object $X\in\lst $ 
the map $ \eta_{X}$ is an $ S$-local equivalence{\rm .} 
\endproclaim

{\it Proof of Theorem~{\rm \ref{NEWMAIN}} assuming the lemma}.
We need to show that if $ X$ is a cofibrant object
in $ \lst$, $ A$ is fibrant in $ \alg$ and 
$ f\colon X\ra A$ is a map in $ \lst$, then 
$ f$ is an $ S$-local equivalence if and only if
its adjoint $f^\flat \colon K_{\T}X\ra A$
is a weak equivalence in $ \alg$ (that is, an
objectwise weak equivalence).
Let $ f$ be a map as above.
We have a commutative diagram
\figin{p16}{900}%
where $ \eta_X$ is an $ S$-local equivalence
by Lemma \ref{UNIT OF ADJUNCTION}. If we assume that 
$ f$ is an $ S$-local equivalence, then so must be
$ f^\flat$. But  $f^\flat$ is a map of 
strict $ \T$-algebras, so  (\ref{LOCAL W.E.})
implies that it is an objectwise weak equivalence. 
Conversely, if $ f^\flat$ is a weak equivalence in
$ \alg$, then it is also an $ S$-local equivalence,
and $ f= f\circ\eta_X$ is an $ S$-local equivalence 
as well. 
\hfill\qed
\vglue12pt 

{\it Proof of Lemma} \ref{UNIT OF ADJUNCTION}.
By (\ref{S-EQUIV. VIA SIMPL. COMPLEXES })
it suffices to show that for every homotopy 
$ \T$-algebra $ Z$ fibrant in $ \stc$ and 
for every cofibrant $X\in\stf$ the map 
$ \eta_X \colon X\ra K_\T X$ induces a weak
equivalence of simplicial function complexes
$$\eta_X^\ast\colon\map(K_\T X,Z)\ra\map(X,Z).$$
The proof is split into a few steps.

1) $X=\coprod_{i=1}^m F_{n_i}$.  
Recall that by (\ref{COR TO K_T OF COPRODUCTS})
we have
$ K_\T(\coprod_{i=1}^m F_{n_i})\simeq F_{\Sigma n_i}$.
For any $ Z\in\stc$, 
$$\map\left(\coprod_{i=1}^m F_{n_i},Z\right)\cong
\prod_{i=1}^m\map(F_{n_i}, Z)\cong 
\prod_{i=1}^m Z\left(T_{n_i}\right)$$
and
$$\map(F_{\Sigma{n_i}}, Z)\cong Z(T_{\Sigma{n_i}}) .$$
Moreover, if $ Z$ is a homotopy $ \T$-algebra 
then $Z(T_{\Sigma{n_i}})\simeq Z(T_1)^{\Sigma{n_i}}
\simeq \prod_{i=1}^m Z(T_{n_i})$. One can check that the 
above weak equivalence is in fact induced by the map 
$\eta_X$ (notice that in this case 
$\eta_X=\kappa$ where $ \kappa $ is the map 
as in (\ref{K_T OF COPRODUCTS})).  

2) $X=\coprod_{i\in I} F_{n_i}$ -- possibly 
infinite disjoint union of free strict $\T$-algebras. 
Let ${\bo P}_I$ denote the category of all 
finite subsets of $ I$ with inclusions of sets
as morphisms. Define a functor 
$$\widetilde{X}\colon {\bo P}_I\ra\lst,\ \ \ \ \ 
\widetilde{X}(A):=\coprod_{i\in A} F_{n_i}.$$
Then $ {\rm colim\, }\widetilde{X}=X$ and 
${\rm colim\, }J_\T K_\T \widetilde{X}=K_\T X$.
The second equality follows from the fact that 
$ J_\T$ commutes with filtered colimits 
(\ref{FILTEREDCOLIMS}) and $ K_\T$ as a left 
adjoint functor preserves all colimits. 
The map $ \eta_X$ is then a colimit of
maps $\eta_{\widetilde X(A)}\colon \widetilde X(A)
\ra K_\T \widetilde X(A)$. For every 
$ A\in{\bo P}_I$ the diagram $ \widetilde{X}(A)$
is of the form   in step 1, so that the map 
induced by $\eta_{\widetilde X(A)} $ on 
simplicial function complexes $\map(-,Z)$ is a weak equivalence. 
It follows that $ {\rm hocolim\, }\eta_{\widetilde X}$
also induces a weak equivalence
$$ {\rm hocolim\, }\eta_{\widetilde X}^\ast
\colon\map({\rm hocolim\, }J_\T K_\T\widetilde X,Z)\,
\mapright{20}{\simeq} \,\map({\rm hocolim\, }\widetilde X,Z).$$
Therefore it is enough to show that 
$ {\rm colim\, }\widetilde{X}\simeq
{\rm hocolim\, }\widetilde{X}$ and
$ {\rm colim\, }J_\T K_\T\widetilde{X}\simeq
{\rm hocolim\, }J_\T K_\T\widetilde{X}$.
By the definition of homotopy colimits 
\cite[19.1.2]{hirschhorn} and since 
the simplicial structure on $ \stc$ is defined 
objectwise,
$$ ({\rm hocolim\, }\widetilde X)(T_n)\cong
{\rm hocolim\, }\widetilde X(T_n).$$
But $ {\bo P}_I$ is a filtered category and by  
\cite[3.5, p.331]{bousfield-kan} ordinary and 
homotopy colimits over filtered categories 
coincide, so that  $ ({\rm hocolim\, }\widetilde X)(T_n)
\simeq({\rm colim\, }\widetilde X)(T_n)$.
Hence $ {\rm hocolim\, }\widetilde X\simeq
{\rm colim\, }\widetilde X$.
Similarly $ {\rm hocolim\, }J_\T K_\T\widetilde X
\simeq {\rm colim\, }J_\T K_\T\widetilde X$.

3) $X=\coprod_{i\in I} F_{n_i}\otimes K_i$ 
($K_i$ - simplicial set). 
Let $ X_\bullet $ denote a simplicial object 
in $ \stc$ such that 
$ X_k := \coprod_{i\in I}\coprod_{\sigma\in (K_i)_k}
F_{n_i}$. Then \cite[VII 3.7]{goerss-jardine}
$ |X_\bullet|\cong X$. Also by (\ref{REALIZATION OF K_T})
we see that $|K_\T X_\bullet|\cong K_\T X$. 
Since for all $ k\geq~0$ the diagram $ X_k$ is of the 
form considered in step 2,  $ \eta_{X_k}$
must be an\break $ S$-local equivalence. Therefore, by 
(\ref{W.E. FOR SIMPLICIAL OBJECTS}) 
and (\ref{4.2 FOR SCC}) the map 
$ \eta_X $ also is an $ S$-local equivalence. 

4) $X=|FU_\bullet Y|$ (see \ref{RESOLUTION}) where 
$ Y$ is any object of $ \stc$.
Since $ FU_k Y$ is of the form in step 3
for any $ k\geq 0$, we can use (\ref{REALIZATION OF K_T}) 
and  (\ref{W.E. FOR SIMPLICIAL OBJECTS}) as above
to show that $ \eta_X$ must be an $ S$-local equivalence.

5) Let $X$ be an arbitrary diagram cofibrant in $ \stf$.
Recall (\ref{AUGMENTATION IS W.E.}) that we have 
a weak equivalence 
$ |\varphi|\colon |FU_\bullet X|\ra X$.
Consider the commutative diagram 
\figin{p18}{900}%
The map $ \eta_{|FU_\bullet X|}$ is an $ S$-local 
equivalence by step 4,
so that  it suffices to show that $ K_{\T}|\varphi|$
is an $ S$-local equivalence.
The functor $ K_{\T}\colon\lst\ra\alg$ is a left 
adjoint in a Quillen pair, and so it preserves all
acyclic cofibrations between cofibrant objects.
Therefore, by K. Brown's lemma 
\cite[9.9]{dwyer-spalinski}
it  preserves all 
weak equivalences between cofibrant objects. 
Since $ |\varphi|$ is a weak equivalence 
in $ \lst$ and both $ X$ and $ |FU_\bullet X|$
are cofibrant we get that $ K_{\T}|\varphi|$
is a  weak equivalence.
\hfill\qed

\end{document}